\documentclass[12pt]{amsart}
\usepackage{amsmath,amssymb,amsthm,amscd}
\newtheorem{formula}{}[section]
\newtheorem{proposition}[formula]{Proposition}
\newtheorem{corollary}[formula]{Corollary}
\newtheorem{lemma}[formula]{Lemma}
\newtheorem{theorem}[formula]{Theorem}
\theoremstyle{definition}
\newtheorem{definition}[formula]{Definition}
\newtheorem{example}[formula]{Example}
\theoremstyle{remark}
\newtheorem*{remark}{Remark}

\def\A{{\mathcal A}}

\begin{document}

\title[Formal vector fields and Buchstaber's conjecture]
{Algebra of formal vector fields on the line and Buchstaber's
conjecture} \subjclass{17B56; 55S30}
\author{Dmitri Millionschikov}
\address{Moscow State University, Department of Mathematics and Mechanics,
Leninskie gory, 119899 Moscow, Russia}
\email{million@mech.math.msu.su}
\date{December 25, 2005}
\keywords{Massey products, graded Lie algebras, formal connection,
Maurer-Cartan equation, representation, cohomology}
\thanks{The research of the author was partially supported by grants
RFBR 05-01-01032 and "Russian Scientific Schools"}
\begin{abstract}
Let $L_1$ denotes the Lie algebra of formal vector fields on the
line which vanish at the origin together with their first
derivatives. $L_1$ is a nilpotent "positive part" of the Witt
(Virasoro) algebra. Buchstaber and Shokurov have shown that the
universal enveloping algebra $U(L_1)$ is isomorphic to the tensor
product $S\otimes {\mathbb R}$, where $S$ is the Landweber-Novikov
algebra  in complex cobordism theory. Goncharova calculated the
cohomology $H^*(L_1)=H^*(U(L_1))$, in particular it follows from
her theorem that $H^*(L_1)$ has trivial multiplicative structure.
Buchstaber conjectured that $H^*(L_1)$ is generated with respect
to non-trivial Massey products by $H^1(L_1)$. Feigin, Fuchs and
Retakh found representation of $H^*(L_1)$ by trivial Massey
products. Later Artelnykh found non-trivial Massey products for a
part of $H^*(L_1)$. In the present article we prove that
$H^*(L_1)$ is generated with respect to non-trivial Massey
products by two elements from $H^1(L_1)$.
\end{abstract}
\date{}

\maketitle

\section*{Introduction}
Buchstaber and Shokurov discovered \cite{BuSho} that the
Landweber-Novikov algebra in the complex cobordism theory tensored
by real numbers $S \otimes {\mathbb R}$ is isomorphic to the
universal enveloping algebra $U(L_1)$ of the Lie algebra $L_1$ of
polynomial vector fields on the real line ${\mathbb R}^1$ with
vanishing non-positive Fourier coefficients. $L_1$ is a maximal
(residually) nilpotent subalgebra of the Witt (Virasoro) algebra.
In that time the algebra $L_1$ attracted a lot of interest
\cite{Fu} and the computation of $H^*(L_1)$ by Goncharova \cite{G}
was one of the most technically complicated results in homology
algebra. Her result allowed Buchstaber and Kholodov to obtain some
deep results in the complex cobordism theory.

It follows from the Goncharova theorem \cite{G} that the
cohomology algebra $H^*(L_1)$ has a trivial multiplication.
Buchstaber conjectured that the algebra $H^*(L_1)$ is generated
with respect to the non-trivail Massey products by its  first
cohomology $H^1(L_1)$.

Feigin, Fuchs and Retakh \cite {FeFuRe} represented the basic
homogeneous cohomology classes from $H^*(L_1)$ as  Massey products
\cite{FeFuRe}. But all the products considered by them are trivial
ones. Twelve years later Artel'nykh \cite{Artel} represented a
part of basic cocycles in $H^*(L_1)$ by means of non-trivial
Massey products, but his brief article contains no proof.

In the present article we prove Buchstaber's conjecture in its
original setting. The main result is the Theorem \ref{main_th}
stating that {\it the cohomology algebra $H^*(L_1)$ is generated
with respect to the non-trivial Massey products by $H^1(L_1)$}.

Although we have strengthened the Feigin-Fuchs-Retakh theorem we
use some important technical tools from \cite{FeFuRe}. One of them
is the free resolution of the trivial $L_1$-module constructed by
means of so-called singular Virasoro vectors (an analogue of
Bernshtein-Gelfand-Gelfand resolution). At the time of writing
\cite{FeFuRe} there were no general formula for singular vectors
in terms of operators $S_{p,q}(t)$ from $U(L_1)$. The formula for
$S_{p,q}(t)$ obtained later by Benoit and Saint-Aubin is another
important ingredient of our proof.

The most important part of our construction is the new graded
so-called thread $L_1$-module $\tilde M$ that we introduce in the
Section \ref{new_module}. The structure and important properties
of  $\tilde M$ is the cornerstone of the proof of the main Theorem
\ref{main_th}.

In the Section \ref{Massey_S} we present May's approach to the
definition of Massey products, his notion of formal connection
developed by Babenko and Taimanov in \cite{BaTa} for Lie algebras,
we introduce also the notion of equivalent Massey products. The
analogy with the classical Maurer-Cartan equation is especially
transparent in the case of Massey products of $1$-dimensional
cohomology classes $\langle \omega_1,\dots,\omega_n\rangle$. The
relation of this special case to the representations theory was
discovered in \cite{FeFuRe}, \cite{Dw}.  Following \cite{FeFuRe}
we consider Massey products $\langle
\omega_1,\dots,\omega_n,\Omega\rangle$, where
$\omega_1,\dots,\omega_n$ are closed $1$-forms and $\Omega$ is a
closed $q$-form.

We stress on non-triviality of our Massey products. It was pointed
out by May to the author that it follows from some general results
\cite{JPM, JPM2} that $H^*(L_1)$ is generated by {\it matrix}
(possibly trivial) Massey products and $H^1(L_1)$. The triple
non-trivial Massey product in cohomology $H^*(L_1)$ was used by
Babenko and Taimanov in their construction of simply connected non
formal symplectic manifolds \cite{BaTa}.

\section{Cohomology of ${\mathbb N}$-graded Lie algebras}

Let $\mathfrak{g}$ be a Lie algebra over ${\mathbb K}$ and
$\rho: \mathfrak{g} \to \mathfrak{gl}(V)$ its linear representation
(or in other words $V$ is a $\mathfrak{g}$-module).
We denote by $C^q(\mathfrak{g},V)$
the space of $q$-linear skew-symmetric mappings of $\mathfrak{g}$ into
$V$. Then one can consider an algebraic complex:

$$
\begin{CD}
V @>{d_0}>>
C^1(\mathfrak{g}, V) @>{d_1}>> C^2(\mathfrak{g}, V) @>{d_2}>>
\dots @>{d_{q{-}1}}>> C^q(\mathfrak{g}, V) @>{d_q}>> \dots
\end{CD}
$$
where the differential $d_q$ is defined by:

\begin{equation}
\begin{split}
(d_q f)(X_1, \dots, X_{q{+}1})=
\sum_{i{=}1}^{q{+}1}(-1)^{i{+}1}
\rho(X_i)(f(X_1, \dots, \hat X_i, \dots, X_{q{+}1}))+\\
+ \sum_{1{\le}i{<}j{\le}q{+}1}(-1)^{i{+}j{-}1}
f([X_i,X_j],X_1, \dots, \hat X_i, \dots, \hat X_j, \dots, X_{q{+}1}).
\end{split}
\end{equation}

The cohomology of the complex $(C^*(\mathfrak{g}, V), d)$ is called
the cohomology of the Lie algebra $\mathfrak{g}$
with coefficients in the representation $\rho: \mathfrak{g} \to V$.

The cohomology of $(C^*(\mathfrak{g}, {\mathbb K}), d)$
($V= {\mathbb K}$ and $\rho: \mathfrak{g} \to {\mathbb K}$ is trivial)
is called the
cohomology with trivial coefficients of the Lie algebra
$\mathfrak{g}$ and is denoted by $H^*(\mathfrak{g})$.

One can remark that
$d_1: C^1(\mathfrak{g}, {\mathbb K}) \rightarrow
C^2(\mathfrak{g}, {\mathbb K})$ of the $(C^*(\mathfrak{g}, {\mathbb K}), d)$
is the dual mapping to the Lie bracket
$[ \, , ]: \Lambda^2 \mathfrak{g} \to \mathfrak{g}$.
Moreover the condition $d^2=0$ is equivalent to the Jacobi identity for $[,]$.

\begin{definition}
A Lie algebra $\mathfrak{g}$ is called $\mathbb N$-graded,
if it is decomposed to the direct sum of subspaces such that
$$\mathfrak{g}=\oplus_{i} \mathfrak{g}_{i}, \; i \in
{\mathbb N},
\quad \quad [\mathfrak{g}_{i}, \mathfrak{g}_{j}]
\subset \mathfrak{g}_{i+j}, \;
\forall \: i, j \in {\mathbb N}.
$$
\end{definition}

\begin{example}
Let us recall that the Witt algebra $W$ is
spanned by differential operators on the real line ${\mathbb R}^1$
with a fixed coordinate $x$
$$
e_i=x^{i+1}\frac{d}{dx}, \; i \in {\mathbb Z}, \quad
[e_i,e_j]= (j-i)e_{i{+}j}, \; \forall \;i,j \in {\mathbb Z}.
$$
We denote by $L_1$ a positive part of the Witt algebra, i.e.
$L_1$ is a subalgebra of $W$ spanned by all $e_i, \; i \ge 1$.

Obviously $W$ is a ${\mathbb Z}$-graded Lie algebra
with one-dimensional homogeneous
components:
$$
W=\oplus_{i\in {\mathbb Z}}W_i, \; W_i=Span(e_i).
$$
Thus $L_1$ is a ${\mathbb N}$-graded Lie algebra.
\end{example}

Let $\mathfrak{g}=\oplus_{\alpha}\mathfrak{g}_{\alpha}$ be a
${\mathbb Z}$-graded Lie algebra and $V=\oplus_{\beta} V_{\beta}$
is a ${\mathbb Z}$-graded $\mathfrak{g}$-module, i.e.,
$\mathfrak{g}_{\alpha}V_{\beta} \subset V_{\alpha {+} \beta}.$
Then the complex $(C^*(\mathfrak{g}, V), d)$ can be equipped with
the ${\mathbb Z}$-grading $C^q(\mathfrak{g},V) = \bigoplus_{\mu}
C^q_{(\mu)}(\mathfrak{g},V)$, where a $V$-valued $q$-form $c$
belongs to $C^q_{(\mu)}(\mathfrak{g},V)$ if and only if for $X_1
\in \mathfrak{g}_{\alpha_1}, \dots,  X_q \in
\mathfrak{g}_{\alpha_q}$ we have
$$c(X_1,\dots,X_q) \in
V_{\alpha_1{+}\alpha_2{+}\dots{+}\alpha_q{+}\mu}.$$

This grading is compatible with the differential  $d$ and
hence we have ${\mathbb Z}$-grading in the cohomology:
$$
H^{q} (\mathfrak{g},V)= \bigoplus_{\mu \in {\mathbb Z}}
H^{q}_{(\mu)} (\mathfrak{g},V).
$$

\begin{remark}
The trivial $\mathfrak{g}$-module ${\mathbb K}$ has only one non-trivial
homogeneous component ${\mathbb K}={\mathbb K}_0$.
\end{remark}

The exterior product in $\Lambda^*(\mathfrak{g})$
induces a structure of a
bigraded algebra in the cohomology
$H^*(\mathfrak{g})$:
$$
H^{q}_{k} (\mathfrak{g}) \wedge
H^{p}_{l} (\mathfrak{g}) \to
H^{q{+}p}_{k+l} (\mathfrak{g}).
$$

Let $\mathfrak{g}= \oplus_{\alpha >0} \mathfrak{g}_{\alpha}$ be a
${\mathbb N}$-graded Lie algebra and $V$ be a
$\mathfrak{g}$-module provided with the following invariant flag
of linear subspaces $V_i, \in {\mathbb Z}$:
$$
V=V_{i_0}\supset V_{i_0+1} \dots \supset V_i \supset V_{i+1}
\supset \dots,
$$
for some $i_0 \in {\mathbb Z}$ and such that $\cap_i V_i=\left\{ 0
\right\}$ and
$$
{\mathfrak g}V_i \subset V_{i+1}, \; i \in {\mathbb Z}.
$$

One can define a decreasing filtration ${\mathcal F}$ of
$(C^*(\mathfrak{g},V),d)$:
$$
{\mathcal F}^0 C^*(\mathfrak{g},V)\supset \dots
\supset {\mathcal F}^q C^*(\mathfrak{g},V)
\supset {\mathcal F}^{q{+}1} C^*(\mathfrak{g},V) \supset \dots
$$
where the subspace ${\mathcal F}^q C^{p{+}q}(\mathfrak{g},V)$ is
spanned by $(p{+}q)$-forms $c$ in $C^{p{+}q}(\mathfrak{g},V)$ such
that
$$
c(X_1,\dots,X_{p{+}q}) \in  V_q, \; \forall X_1,\dots,X_{p{+}q}
\in \mathfrak{g}.
$$

The filtration ${\mathcal F}$ is compatible with $d$.

Let us consider the corresponding spectral sequence $E_r^{p,q}$:
\begin{proposition}
\label{spectral_seq}
$E_1^{p,q}=\left(V_q/V_{q+1} \right)\otimes
H^{p{+}q}({\mathfrak g})$.
\end{proposition}
\begin{proof}
We have the following natural isomorphisms:
\begin{equation}
\begin{split}
C^{p{+}q}({\mathfrak g},V)
= V \otimes \Lambda^{p{+}q}({\mathfrak g}^*)\\
E_0^{p,q}={\mathcal F}^q C^{p{+}q}(\mathfrak{g},V)/ {\mathcal
F}^{q{+}1} C^{p{+}q}(\mathfrak{g},V) =(V_q/ V_{q+1}) \otimes
\Lambda^{p{+}q}({\mathfrak g}^*).
\end{split}
\end{equation}

Now the proof follows from the formula for the
$d_0^{p,q}: E_0^{p,q} \to E_0^{p{+}1,q}$:
$$d_0(v \otimes f)=v\otimes df,$$
where $v \in V, f \in \Lambda^{p{+}q}({\mathfrak g}^*)$
and $df$ is the standard
differential of the cochain complex of $\mathfrak{g}$ with trivial
coefficients.
\end{proof}

\begin{theorem}[Goncharova,\cite{G}]
The Betti numbers ${\rm dim} H^q(L_1)=2$, for every
$q \ge 1$, more precisely
$$ {\rm dim} H_k^q(L_1)=
\left\{\begin{array}{r}
   1, \hspace{0.6em}{\rm if}~k=\frac{3q^2 \pm q}{2}, \\
   0,\hspace{1.76em}{\rm otherwise.}\\
   \end{array} \right . \hspace{3.3em} $$
\end{theorem}
We will denote in the sequel by $g^q_{\pm}$ the generators of the
spaces $H^q_{\frac{3q^2\pm q}{2}}(L_1)$. The numbers $\frac{3q^2
\pm q}{2}$ that we will denote sometimes in the sequel by
$e_{\pm}(q)$ are so called Euler pentagonal numbers. A sum of two
arbitrary pentagonal numbers is not a pentagonal number, hence the
algebra $H^*(L_1)$ has a trivial multiplication.

\begin{example}

1) $H^1(L_1)$ is generated
by $g^1_-=[e^1]$ and $g^1_+=[e^2]$;

2) the basis of $H^2(L_1)$ consists of two
classes $g^2_-=[e^1 {\wedge} e^4]$ and
$g^2_+=[e^2 {\wedge} e^5-3e^3 {\wedge} e^4]$ of weights $5$ and $7$
respectively.
\end{example}
\begin{remark}
It is very difficult (almost impossible) to understand all the
details of Goncharova's proof. There is another one proof by
Weinstein \cite{W1, w2}, there is a mistake in \cite{W1} corrected
later in \cite{W2} in the construction of the spectral sequence,
but the case $L_1$ can be treated apart from the general case as
it was shown in \cite{GeFeFu} (corrected later in \cite{W1}).
\end{remark}

\section{Massey products in cohomology.}
In this section we follow \cite{JPM} and \cite{BaTa} presenting the
definitions of Massey products.
\label{Massey_S}
Let $\A=\oplus_{l \ge 0}\A^l$
be a differential graded algebra over a field ${\mathbb K}$.
It means that the following operations are defined:
an associative multiplication
$$
\wedge:\A^l\times\A^m\to\A^{l+m},\; l,m\geq 0, \; l,n \in {\mathbb Z}.
$$
such that $a\wedge b=(-1)^{lm}\,b\wedge a$ for $a\in\A^l$, $b\in\A^m$,
and a differential $d, \; d^2=0$
$$
d:\A^l\to\A^{l+1},\ \ l\geq 0,
$$
satisfying the Leibniz rule
$d\,(a\wedge b)=d\,a\wedge b+(-1)^l a\wedge d\,b$ for $a\in\A^l$.
\begin{example}
$\A=\Lambda^*(\mathfrak{g})$ is the cochain complex of a Lie algebra.
\end{example}

For a given differential graded algebra  $(\A,d)$ we denote by
$LT_n(\A)$ a space of all lower triangular $(n+1)\times
(n+1)$-matrices with entries from $\A$, vanishing at the main
diagonal. $LT_n(\A)$ has a structure of a differential algebra
with a standard matrix multiplication, where matrix entries are
multiplying as elements of $\A$. A differential $d$ on $LT_n(\A)$
is defined by
\begin{equation}
d\,A=(d\,a_{ij})_{1\leq i,j\leq n{+}1}.
\label{1.2.1}
\end{equation}

An involution $a\to\bar{a}=(-1)^{k+1}a, a \in A^k$ of $\A$ can be
extended to an involution of $LT_n(\A)$ as
$\bar{A}=(\bar{a}_{ij})_{1\leq i,j\leq n{+}1}$. It satisfies the
following properties:
$$
\overline{\bar{A}}=A, \quad \overline{AB}=-\bar{A}\bar{B}, \quad
\overline{d\,A}=-d\,\bar{A}.
$$
Also we have the generalized Leibniz rule for the differential
(\ref{1.2.1})
$$
d\,(AB)=(d\,A)B-\bar{A}(d\,B).
$$

The algebra $LT_n(\A)$ has a two-sided center $I_n(\A)$ of
matrices
$$
\left(\begin{array}{cccc}
0 & \dots & 0 & 0 \\
0 & \dots & 0 & 0 \\
& \dots& & \\
\tau & \dots & 0 & 0
\end{array}\right), \quad \tau \in \A
$$

\begin{definition}[\cite{BaTa}]
A matrix $A \in LT_n(\A)$ is called the matrix of a formal
connection if it satisfies the Maurer-Cartan equation
\begin{equation}
\mu(A)=d\,A-\bar{A}\cdot A \in I_n(\A).
\label{star}
\end{equation}
\end{definition}
\begin{proposition}[\cite{BaTa}]
Let $A$ be the matrix of a formal connection,
then the entry $\tau \in \A$ of the matrix $\mu(A) \in I_n(\A)$
in the definition (\ref{star}) is closed.
\end{proposition}
\begin{proof}
We have the following generalized Bianci identity
for the Maurer-Cartan operator
$\mu(A)=d\,A - \bar{A} \cdot A$ ($A$ is an arbitrary matrix):
$$
d\,\mu(A)=\overline{\mu(A)}\cdot A+A\cdot\mu(A).
$$
Indeed it's easy to verify the following equalities:
$$
d\,\mu(A)=-d\,(\bar{A}\cdot A)= -d\,\bar{A}\cdot A+A\cdot
d\,A=\overline{d\,A}\cdot A+A\cdot d\,A=
$$
$$
=\overline{(\mu(A)+\bar{A}\cdot A)}\cdot A+A(\mu(A)+\bar{A}\cdot
A)=
$$
$$
= \overline{\mu(A)}\cdot A- A\cdot\bar{A}\cdot
A+A\cdot\mu{A}+A\cdot\bar{A}\cdot A=
$$
$$
=\overline{\mu(A)}\cdot A+A\cdot\mu(A).
$$
\end{proof}
Now let $A$ be the matrix of a formal connection, then
the matrix $\mu(A)$ belongs to the center $I_n(\A)$ and hence $d\mu(A)=0$.
One can think of $\mu(A)$ as the curvature
matrix of a formal connection $A$.

Let $A$ be an lower triangular matrix from $LT_n(\A)$. One can
rewrite it in the following notation:
$$
A=\left(\begin{array}{cccccc}
  0      & 0 & \dots & 0 & 0   & 0   \\
  a(n,n) & 0      & \dots & 0 & 0   & 0   \\
  a(n-1,n)  & a(n-1,n-1)  & \dots  & 0 & 0      & 0    \\
  \dots      & \dots      & \dots      & \dots &  \dots & \dots  \\
  a(2,n)      & a(2,n-1)      & \dots      & a(2,2) & 0          & 0   \\
  a(1,n)      & a(1,n-1)      & \dots      & a(1,2) & a(1,1)          & 0
  \end{array}\right).
$$

\begin{proposition}
A matrix $A \in LT_n(\A)$ is the matrix of a formal connection if
and only if the following conditions on its entries hold on
\begin{equation}
\label{def_syst}
\begin{split}
a(i,i)=a_i \in \A^{p_i}, \quad i=1,\dots,n;\\
a(i,j)\in\A^{p(i,j)+1}, \quad p(i,j)=\sum^j_{r=i}(p_r-1);\\
d\,a(i,j)=\sum_{r=i}^{j-1}\bar{a}(i,r)\cdot a(r+1,j),\;\;
(i,j)\neq(1,n).
\end{split}
\end{equation}
\end{proposition}
\begin{proof}
The system (\ref{def_syst}) is  just the Maurer-Cartan equation
rewritten in terms of the entries of the matrix $A$ and it
is a part of the classical definition \cite{K} of
the defining system for a Massey product.
\end{proof}
\begin{definition}[\cite{K}]
A collection of elements,
$A=(a(i,j))$, for $1\leq i\leq j\leq n$ and $(i,j)\neq(1,n)$
is said to be a defining  system for the product
$\langle a_1,\dots,a_n\rangle$ if it satisfies (\ref{def_syst}).

In this situation
the $(p(1,n)+2)$-dimensional cocycle
$$
c(A)=\sum_{r=1}^{n-1}\bar{a}(1,r)a(r+1,n)
$$
is called the related cocycle of the defining system $A$.
\end{definition}
\begin{remark}
We saw that the notion of the defining system is equivalent to
the notion of the formal connection. However one has to remark
that an entry $a(1,n)$ of the matrix $A$ of a formal connection
does not belong to the corresponding defining system $A$,
it can be taken as an arbitrary element from $\A$.
In this event the only one (possible) nonzero entry
$\tau$ of the Maurer-Cartan matrix $\mu(A)$ is equal to $-c(A)+da(1,n)$.
\end{remark}

\begin{definition}[\cite{K}]
The $n$-fold product $\langle a_1,\dots,a_n\rangle$ is defined if
there is at least one defining system for it (a formal connection
$A$ with entries $a_n,\dots,a_1$ at the second diagonal). If it is
defined, then $\langle a_1,\dots,a_n\rangle$ consists of all
cohomology classes $\alpha \in H^{p(1,n)+2}(\A)$ for which there
exists a defining system $A$ (a formal connection $A$) such that
$c(A)$ ($-\tau$ respectively) represents $\alpha$.
\end{definition}

\begin{theorem}[\cite{K},\cite{BaTa}]
The operation $\langle a_1,\dots,a_n\rangle$ depends only on the
cohomology classes of the elements $a_1,\dots,a_n$.
\end{theorem}
\begin{proof}
A changing of an arbitrary entry $a_{ij}, \:j > i$
of the matrix $A$ of a formal connection to $a_{ij}+db$ leads
to a replacement of $A$ by
$$
A'=A+db\cdot E_{ij}+A\cdot b\cdot E_{ij}-\bar b\cdot E_{ij}\cdot A,
$$
where $E_{ij}$ is a scalar matrix wich has $1$
on $(i,j)$-th place and zeroes on all others.
For the corresponding Maurer-Cartan matrix we will have
$$
\mu(A')=\mu(A)+d((A\cdot b \cdot E_{ij}-\bar b \cdot E_{ij}\cdot A)\cap I_n).
$$
\end{proof}
\begin{definition}[\cite{K}]
A set of closed elements $a_i, i=1,\dots,n$ from $\A$ representing
some cohomology classes ${\alpha}_i \in H^{p_i}(\A), i=1,\dots,n$
is said to be a defining system for the Massey $n$-fold product
$\langle {\alpha}_1,\dots,{\alpha}_n\rangle$ if it is one for
$\langle a_1,\dots,a_n\rangle$. The Massey $n$-fold product
$\langle {\alpha}_1,\dots,{\alpha}_n\rangle$ is defined if
$\langle a_1,\dots,a_n\rangle$ is defined, in which case $\langle
{\alpha}_1,\dots,{\alpha}_n\rangle=\langle a_1,\dots,a_n\rangle$
as subsets in $H^{p(1,n)+2}(\A)$.
\end{definition}

\begin{example}
For $n=2$ the matrix $A$ of a formal connection has a form
$A=\left(\begin{array}{ccc}
0 & 0 & 0\\
b & 0 & 0\\
c & a & 0
\end{array}\right)$ and the matrix Maurer-Cartan equation is equivalent
to two equations $db=0$ and $da=0$. Evidently $\langle \alpha,
\beta \rangle = \langle a, b \rangle =\bar \alpha \cdot \beta$.
\end{example}

\begin{example}[Triple Massey products]
Let $\alpha$, $\beta$, and  $\gamma$ be the cohomology classes of closed
elements $a \in\A^p$, $b\in\A^q$, and  $c\in\A^r$.
The Maurer-Cartan equation for
$$
A=\left(\begin{array}{cccc}
0 & 0 & 0 & 0 \\
c & 0 & 0 & 0 \\
g & b & 0 & 0 \\
h & f & a & 0
\end{array}\right).
$$
is equivalent to
\begin{equation}
d\,f=(-1)^{p+1}\,a\wedge b,\ \ d\,g=(-1)^{q+1}\,b\wedge c.
\label{ast}
\end{equation}
Hence the triple Massey product
$\langle\alpha,\beta,\gamma\rangle$
is defined if and only if
$$
\alpha\cdot \beta=\beta \cdot \gamma=0\ \ \mbox{in}\ \ H^{\ast}(\A).
$$
If these conditions are satisfied then
the Massey product $\langle \alpha, \beta, \gamma \rangle$ is defined
as a subset in $H^{p{+}q{+}{r}{-}1}(\A)$ of the following form
$$
\langle \alpha, \beta, \gamma \rangle=\left\{ [(-1)^{p+1} a\wedge
g+(-1)^{p+q} f\wedge c]\right\}.
$$
Since $f$ and  $g$ are defined by (\ref{ast}) up to closed elements
from $\A^{p+q-1}$  and $\A^{q+r-1}$ respectively, the triple Massey product
$\langle \alpha,\beta,\gamma\rangle$ is an affine subspace of
$H^{p{+}q{+}{r}{-}1}(\A)$ parallel to
$\alpha \cdot H^{q+r-1}(\A)+ H^{p+q-1}(\A)\cdot \gamma$.
\end{example}
\begin{remark}
We defined Massey products as the multi-valued operations in
general. More often in the literature the triple Massey product is
defined as a quotient-space $\langle \alpha, \beta, \gamma
\rangle/ (\alpha \cdot H^{q+r-1}(\A)+ H^{p+q-1}(\A)\cdot \gamma)$
and it is single-valued in this case (see \cite{Fu}).
\end{remark}

\begin{definition}
Let an  $n$-fold Massey product
$\langle {\alpha}_1,\dots,{\alpha}_n\rangle$
be defined. It is called trivial if it contains the trivial cohomology class:
$0\in\langle {\alpha}_1,\dots,{\alpha}_n\rangle$.
\end{definition}

\begin{proposition}
\label{triviality}
Let a Massey product
$\langle {\alpha}_1,\dots,{\alpha}_n\rangle$
is defined. Then all Massey products
$\langle {\alpha}_l,\dots,{\alpha}_q\rangle, 1\le l < q \le n, q-l<n-1$
are defined and trivial.
\end{proposition}
\begin{proof}
It follows from (\ref{def_syst}).
\end{proof}
\begin{remark}
The triviality of all Massey products
$\langle {\alpha}_l,\dots,{\alpha}_q\rangle, 1\le l < q \le n, q-l<n-1$
is only a necessary condition
for a Massey product $\langle {\alpha}_1,\dots,{\alpha}_n\rangle$
to be defined. It is sufficient only in the case $n=3$.
\end{remark}

Let us denote by $GLT_n({\mathbb K})$ a group of non-degenerate
lower triangular $(n{+}1,n{+}1)$-matrices of the form:
$$
C=\left(\begin{array}{ccccc}
  a_{1,1}      & 0 & \dots & 0   & 0   \\
  a_{1,2}      & a_{2,2}      &  \dots & 0  & 0   \\
 \dots   &   \dots     &       \dots &  \dots          & \dots     \\
  a_{1,n}      & a_{2,n}            & \dots & a_{n,n}& 0   \\
  a_{1,n{+}1}      & a_{2,n{+}1}            & \dots & a_{n,n{+}1}          & a_{n{+}1,n{+}1}
  \end{array}\right).
$$
\begin{proposition}
Let $A \in LT_n({\A})$ be the matrix of a formal connection and
$C$ an arbitrary matrix from $GLT_n({\mathbb K})$. Then the matrix
$C^{-1}AC\in LT_n({\A})$ and satisfies the Maurer-Cartan equation,
i.e. is again the matrix of a formal connection.
\end{proposition}
\begin{proof}
$$
d(C^{-1}AC)-\bar C^{-1}\bar A\bar C \wedge C^{-1}AC=
C^{-1}\left(dA-\bar A\wedge A\right)C=0.
$$
\end{proof}
\begin{example}
Let $A \in LT_n({\A})$ be the matrix of a formal connection
(defining system) for a Massey product $\langle \alpha_1,\dots,
\alpha_n\rangle$. Then a matrix $C^{-1}AC$ with
$$
C=\left(\begin{array}{ccccc}
  x_1{\dots}x_{n{-}1}x_n      & 0 & \dots & 0   & 0   \\
  0      & x_1{\dots}x_{n{-}1}      &  \dots & 0  & 0   \\
 \dots   &        &       \dots &            & \dots     \\
  0      & 0            & \dots & x_1 & 0   \\
  0      & 0            & \dots & 0          & 1
  \end{array}\right)
$$
is a defining system for $\langle x_1\alpha_1,\dots,
x_n\alpha_n\rangle= x_1\dots x_n\langle \alpha_1,\dots,
\alpha_n\rangle$.
\end{example}

\begin{definition}
\label{A_equiv} Two matrices $A$ and $A'$ of formal connections
from  $LT_n({\A})$ are equivalent if there exists a matrix $C \in
GL(n{+}1,{\mathbb K})$ such that
$$
A'=C^{-1}AC.
$$
\end{definition}

\begin{example}
Triple products $\langle\alpha,\beta,\gamma\rangle$ and
$\langle x \alpha,y \beta,z \gamma\rangle$, where
$x,y,z \ne 0$, are equivalent with
$$
C=\left(\begin{array}{cccc}
  xyz      & 0 &  0   & 0   \\
  0      & xy      &   0  & 0   \\
  0      & 0            &  x & 0   \\
  0      & 0            &  0          & 1
  \end{array}\right)
$$
and
$$
\langle x \alpha,y \beta,z \gamma\rangle=
xyz\langle \alpha, \beta, \gamma \rangle, \quad x,y,z \in{\mathbb K}.
$$
\end{example}
\begin{remark}
Following the original Massey work \cite{Mass}
some higher order cohomological operations
that we call now Massey products were introduced
in the 60s in \cite{K} and \cite{JPM}.
The relation between Massey products and the
Maurer-Cartan equation was first noticed by May \cite{JPM} and this analogy
was not developed untill \cite{BaTa}.

In the present article we deal only with Massey products of
non-trivial cohomology classes. It is possible to take some of
them trivial, but in this situation is more natural to work with
so-called matrix Massey products that were first introduced by May
\cite{JPM}. This approach was also developed in \cite{BaTa}. We
will not treat this case in the sequel.
\end{remark}

\section{Formal connections and representations}
\label{g-modules} Let $LT_n({\mathbb K})$ be a Lie algebra of
lower triangular $(n+1,n+1)$-matrices over a field $\mathbb K$ of
zero characteristic and $\rho: \mathfrak{g} \to LT_n({\mathbb K})$
be a representation of a Lie algebra $\mathfrak{g}$.

\begin{example}
We take $n=1$ and consider a linear map
$$
\rho: x \in \mathfrak{g} \to
\left(\begin{array}{cc}
0 & 0 \\
a(x)& 0
\end{array}\right).$$
It is evident that
$\rho$ is a Lie algebra homomorphism if and only if
the linear form $a \in {\mathfrak{g}}^*$ is closed
$$
da(x,y)=a([x,y])=a(x)a(y)-a(y)a(x)=0, \; \forall x,y \in \mathfrak{g}.
$$

In other words the matrix
$A=\left(\begin{array}{cc}
0 & 0 \\
a& 0
\end{array}\right)$
satisfies the "strong" Maurer-Cartan equation
$dA-\bar A \wedge A =0$.
\end{example}
\begin{remark}
We recall that we defined in the Section \ref{Massey_S} the
involution of a graded $\A$ as $\bar a=(-1)^{k+1} a, a \in \A^k$.
Thus for a matrix $A$ with entries from ${\mathfrak{g}}^*$ we have
$\bar A =A$. One has to remark that $\bar a$ differs by the sign
from the definition of $\bar a$ in \cite{K}, however in
\cite{JPM2} one meets our sign rule.
\end{remark}

\begin{proposition}
A matrix $A$ with entries from ${\mathfrak{g}}^*$ defines
a representation $\rho: \mathfrak{g} \to T_n({\mathbb K})$
if and only if $A$ satisfies the strong Maurer-Cartan equation
$$dA-\bar A \wedge A=0.$$
\end{proposition}
\begin{proof}
$$
(dA-\bar A \wedge A)(x,y)=A([x,y])-\left[ A(x), A(y)\right],\;\forall x,y \in
\mathfrak{g}.
$$
\end{proof}
\begin{example}
For $n=2$ the matrix $A$ of a representation
$\rho$ has a form
$A=\left(\begin{array}{ccc}
0 & 0 & 0\\
b & 0 & 0\\
c& a& 0
\end{array}\right),$ where $a,b,c \in {\mathfrak{g}}^*$
and the strong Maurer-Cartan equation is equivalent
to the following equations on entries $a,b,c$:
$$da=0, \quad db=0, \quad dc= a \wedge b.$$
\end{example}

The Lie algebra $LT_n({\mathbb K})$ has a one-dimensional center
$I_n({\mathbb K})$ spanned by the matrix
$$
\left(\begin{array}{cccc}
0 & \dots & 0 & 0 \\
0 & \dots & 0 & 0 \\
& \dots& & \\
1 & \dots & 0 & 0
\end{array}\right).
$$
One can consider an one-dimensional central extension
$$
\begin{CD}
0 @>>> \mathbb K \cong I_n({\mathbb K}) @>>> LT_n(\mathbb K)
@>{\pi}>> \tilde LT_n(\mathbb K) @>>> 0.
\end{CD}
$$

\begin{proposition}[\cite{FeFuRe}, \cite{Dw}]
Fixing a Lie algebra homomorphism $\tilde \varphi : \mathfrak{g}
\to \tilde LT_n({\mathbb K})$ is equivalent to fixing a defining
system $A$ with elements from
$\mathfrak{g}^*=\Lambda^1(\mathfrak{g})$. The related cocycle
$c(A)$ is cohomologious to zero if and only if $\tilde \varphi$
can be lifted to a homomorphism $\varphi : \mathfrak{g} \to
LT_n({\mathbb K})$, $\tilde \varphi = \pi \varphi$.
\end{proposition}

Let us consider a Massey product of the form
$$
\langle \omega_1,\dots,\omega_n, \Omega \rangle, \; \omega_i \in
H^1({\mathfrak g}), i=1,\dots,n, \; \Omega \in H^{p}({\mathfrak
g}).
$$
If it is defined then $n$-fold product $\langle
\omega_1,\dots,\omega_n\rangle$ is trivial and the existence of
the homomorphism $\varphi : \mathfrak{g} \to LT_n({\mathbb K})$
means that there is a $(n+1)$-dimensional ${\mathfrak g}$-module
$V$ with a basis $f_1,\dots,f_{n+1}$, such that
$$
{\mathfrak g}f_j \in {\rm Span}(f_{j+1},\dots,f_{n+1}),
j=1,\dots,n; \;{\mathfrak g}f_{n+1}=0.
$$
One can consider the spectral sequence $E_r^{p,q}$ from the
Proposition \ref{spectral_seq} converging to the cohomology
$H^*({\mathfrak g}, V)$ of ${\mathfrak g}$ with coefficients in
$V$.
\begin{theorem}
\label{Massey+Speq} Let ${\mathfrak g}$ be a Lie algebra and a
$(n+1)$-fold Massey product
$$
\langle \omega_1,\dots,\omega_n, \Omega \rangle, \; \omega_i \in
H^1({\mathfrak g}), i=1,\dots,n, \; \Omega \in H^p({\mathfrak g}),
$$
be defined. Let also $A$ be a corresponding formal connection.
Then it exists the $(n+1)$-dimensional ${\mathfrak g}$-module $V$
and the spectral sequence $E_r^{p,q}$ converging to the cohomology
$H^*({\mathfrak g}, V)$ such that
\begin{equation}
\begin{split}
f_1\otimes \Omega \in E_1^{p-1,1}, \\ d_1(f_1\otimes \Omega)=\dots=d_{n-1}(f_1\otimes \Omega)=0, \\
d_n(f_1\otimes \Omega)= f_{n+1}\otimes [c(A)],
\end{split}
\end{equation}
where $d_i: E_i^{p-i,i} \to E_i^{p-i+2,i+1}$ is the $i$-th
differential of the spectral sequence $E_r^{p,q}$ and $c(A)$ is
the cocycle of the formal connection $A$.
\end{theorem}
\begin{proof}
The proof is almost evident, one has to follow only the
definitions. Namely we denote the entries of the matrix $A$ of the
formal connection by
$$
A=\left(\begin{array}{ccccccc}
  0      & 0 & 0 & \dots & 0   & 0 &0  \\
  \Omega & 0 & 0 & \dots & 0   & 0  &0 \\
  \Omega_1  & \omega_n  & 0 & \dots & 0      & 0 &0   \\
  \Omega_2      & a(n{-}1,n)      & \omega_{n{-}1}      & \dots & 0          & 0&0   \\
  \dots      & \dots      & \dots      & \dots &  \dots & \dots & \dots \\
  \Omega_{n{-}1} & a(2,n)      & a(2,n{-}1)      & \dots      & \omega_2 & 0          & 0  \\
  \Omega_n & a(1,n)      & a(1,n{-}1)      & \dots      & a(1,2) & \omega_1          &
  0
  \end{array}\right).
$$
The Maurer-Cartan equation gives us the following system of
equations on the elements $\Omega_i$ standing at the first column
of $A$:
\begin{equation}
\label{def_syst2}
\begin{split}
d\Omega_1=\bar {\omega_n} \wedge \Omega=\omega_n \wedge \Omega,\\
d\Omega_2= \omega_{n-1} \wedge \Omega_1+a(n{-}1,n) \wedge
\Omega,\\
\quad \dots \quad \quad \dots \quad \quad \dots \quad \quad
\dots\quad \quad \dots\quad \quad \dots\\
d\Omega_{n-1}= \omega_2 {\wedge} \Omega_{n-1}+\dots+a(2,n{-}2)
{\wedge} \Omega_2+a(2,n{-}1) {\wedge} \Omega_1+a(2,n) {\wedge}
\Omega.
\end{split}
\end{equation}
We recall that $\Omega_n$ can be an arbitrary form. The
corresponding cocycle $c(A)$ will be equal to
$$
c(A)=\sum_{i=0}^{n-1} a(1,n{-}i) \wedge \Omega_i, \quad
\Omega_0=\Omega, \; a(1,1)=\omega_1.
$$
On the another hand in the cochain complex $C^*({\mathfrak g},V)$
we have
\begin{equation}
\begin{split}
{\bf d}f(X)=[X,f], \; f \in C^0({\mathfrak g},V)=V,\; X \in V;\\
{\bf d}\left(v\otimes \Omega \right)={\bf d}v \wedge
\Omega-v\otimes d\omega, \; v \in V, \; \Omega \in
\Lambda^*({\mathfrak g}),
\end{split}
\end{equation}
where ${\bf d}$ denotes the differential in the complex
$C^*({\mathfrak g},V)$ and the standard $d$ is the differential of
the complex $C^*({\mathfrak g})$ with trivial coefficients. Hence
${\bf d}f_{n+1}=0$ and for $i=1,\dots,n$, we have
$$
{\bf d}f_i=f_{i+1}\otimes \omega_{n-i+1}+f_{i+2}\otimes
a(n{-}i,n-i+1)+\dots+f_{n+1}\otimes a(1,n-i+1).
$$
And
\begin{equation}
\begin{split}
{\bf d}\left(f_1\otimes \Omega\right)=f_2\otimes \omega_n{\wedge
\Omega}+\dots,\\
{\bf d}\left(f_1\otimes \Omega+f_2 \otimes \Omega_1
\right)=f_3\otimes \left(\omega_{n-1} \wedge \Omega_1+a(n{-}1,n)
\wedge \Omega\right)+\dots,\\
\dots\quad\quad\dots\quad\quad\dots\quad\quad\dots\quad\quad\\
{\bf d}\left(\sum_{i=1}^{n}f_i \otimes \Omega_{i-1}
\right)=f_{n+1}\otimes \left( \sum_{i=0}^{n-1} a(1,n{-}i) \wedge
\Omega_i\right)=f_{n+1}\otimes c(A).
\end{split}
\end{equation}
The end of the proof.
\end{proof}

\section{The Feigin-Fuchs-Retakh theorem}

We recall that the algebra $H^*(L_1)$ has a trivial
multiplication. Buchstaber conjectured that the algebra $H^*(L_1)$
is generated with respect to the Massey products by $H^1(L_1)$,
moreover all corresponding Massey products can be chosen non
trivial. The weak version of Buchstaber's conjecture was proved by
Feigin, Fuchs and Retakh \cite{FeFuRe}.

\begin{theorem}[Feigin, Fuchs, Retakh \cite{FeFuRe}]
\label{FeigFu} Let $g^k_{\pm}$ be a non trivial cocycle in
$H^k_{\frac{1}{2}(3q^2\pm q)}(L_1)$. For any $k \ge 2$ we have
$$
g^k_- \in \langle \underbrace{e^1,\dots,e^1}_{2k-1},g^{k{-}1}_+,
\rangle, \quad g^k_+ \in \langle
\underbrace{e^1,\dots,e^1}_{3k-1},g^{k{-}1}_+ \rangle.
$$
\end{theorem}

Feigin, Fuchs and Retakh proposed the following formal connection:
$$
\label{FeFuReCon} A=\left(\begin{array}{ccccccc}
  0      & 0 & 0 & \dots & 0   & 0 &0  \\
  g^{k{-}1}_+ & 0 & 0 & \ddots & \ddots   & 0  &0 \\
  \Omega_1  & e^1  & 0 & \ddots & \ddots      & \ddots &0   \\
  \Omega_2      & \alpha e^2      & e^1      & \ddots & \ddots          & \vdots &\vdots   \\
  \vdots      & 0      & \alpha e^2      & \ddots &  \ddots & 0 & \vdots \\
  \Omega_{n{-}1} & \vdots      & \ddots      & \ddots      & e^1 & 0          & 0  \\
  * & 0      & \ldots      & 0      & \alpha e^2 & e^1  &  0
  \end{array}\right),
$$
with homogeneous forms $\Omega_i \in
\Lambda^{k+i-1}_{\frac{1}{2}(3(k{-}1)^2{+}(k{-}1))+i}(L_1)$ and
parameter $\alpha \in {\mathbb K}$. The system (\ref{def_syst2})
on $\Omega_i, i=1\dots,n{-}1$, is solvable for all values of
parameter $\alpha$ if $n \le 2k$. The corresponding cocycle $c(A)
\in
\Lambda^{k+1}_{\frac{1}{2}(3(k{-}1)^2{+}(k{-}1)){+}n{-}1}(L_1)$
can be non-trivial only if $n=2k$ or $n=3k$ that corresponds to
$H^k_{\frac{1}{2}(3k^2\pm k)}(L_1)$. According to the Theorem
\ref{Massey+Speq} the triviality of the cocycle $c(A)$ is
equivalent to the triviality of the differential $d_n$ of the
spectral sequence $E_r^{p,q}$ converging to $H^*(L_1,V)$.

By methods that we are going to discuss in next sections Feigin,
Fuchs and Retakh established that

1) $d_{2k}$ is trivial if and only if $\alpha \in
\{\frac{1}{6},\frac{1}{24},\dots,\frac{1}{6(k{-}1)^2}\}$;

2) $d_{3k}$ is defined and trivial if and only if:

a) $\alpha \in
\{\frac{1}{6},\frac{1}{54},\dots,\frac{1}{6(k{-}3)^2},\frac{1}{6(k{-}1)^2}\}$
in the case of even $k$;

b) $\alpha \in
\{\frac{1}{24},\frac{1}{96},\dots,\frac{1}{6(k{-}3)^2},\frac{1}{6(k{-}1)^2}\}$
if $k$ is odd.
\begin{corollary}
All Massey products from the Theorem \ref{FeigFu} are trivial.
\end{corollary}

A few words about the proof of the Theorem \ref{FeigFu}. The main
technical problem is obvious: one have to write out explicit
formulas for the forms $\Omega_i$ from the formal connection $A$.
It easy to see that cocycles $g^2_-=e^2{\wedge}e^3$ and
$g^2_+=e^2{\wedge}e^5-e^3{\wedge}^4$ span the homogeneous
subspaces $H^2_5(L_1)$ and $H^2_7(L_1)$ respectively. But it is
still an open question how to write out explicit formulas for all
Goncharova's cocycles $g^k_{\pm}$ in terms of exterior forms from
$\Lambda^*(L_1)$. Fuchs, Feigin and Retakh proposed \cite{FeFuRe}
an elegant and effective way how to present an alternative proof
of non-triviality of differentials (for the special choice of the
parameter $\alpha$) in the spectral sequence $E^{p,q}_r$ without
providing explicit formulas for $\Omega_i$.

In $2000$ Buchstaber's PhD-student Artel'nykh attacked
Buchstaber's conjecture in its original setting. In particular he
claimed the following

\begin{theorem}[Artel'nykh \cite{Artel}]
\label{Th_Artel}
There are non-trivial Massey products
$$
g^k_- \in \langle \underbrace{e^2,\dots,e^2}_{k-1},
g^{k{-}1}_+,e^1 \rangle, k \ge 2, \quad g^{2l{+}1}_+ \in \langle
\underbrace{e^2,\dots,e^2}_{3l+1}, g^{2l}_+ \rangle, l \ge 1.
$$
\end{theorem}
One can see that Artel'nykh have not found non-trivial Massey
products for cohomology classes $g^{2l}_+$. On the another hand
Artel'nykh's brief article \cite{Artel} does not contain any
proof.
\begin{remark}
The inclusion $g^{2l{+}1}_+ \in \langle
\underbrace{e^2,\dots,e^2}_{3l+1}, g^{2l}_+ \rangle, l \ge 1,$ can
be found in \cite{FeFuRe}. The question is the corresponding
Massey product non-trivial or not have not been
discussed.
\end{remark}

\section{Verma modules and Virasoro singular vectors}
\begin{definition}
The Virasoro algebra Vir is an infinite dimensional Lie algebra
defined by its basis $\{z, e_i, i \in {\mathbb Z}\}$ and the
structure relations:
$$
[e_i,z]=0, \; \forall i \in {\mathbb Z}, \quad
[e_i,e_j]=(j-i)e_{i+j}+\frac{j^3-j}{12}\delta_{-i,j}z.
$$
\end{definition}
Vir is ${\mathbb Z}$-graded Lie algebra where $z, e_0$ have
gradings equal to zero, and a generator $e_i$ has grading equal to
$i$. Vir is the one-dimensional central extension of the Witt
algebra (the one-dimensional centre of Vir is spanned by $z$).

\begin{definition}
A Virasoro Verma module $V(h,c)$ is a free module over the
subalgebra $L_1 \subset {\rm Vir}$ generated by some vector $v$
such that
$$
zv=cv, \; e_0v=hv,\quad e_iv=0,\quad {\rm if } \; i <0,
$$
where $c,h \in {\mathbb C}$.
\end{definition}
A Verma module $V(h,c)$ is ${\mathbb N}$-graded module:
$$
V(h,c)=\bigoplus_{n=0}^{+\infty} V_n(h,c), \quad V_n(h,c)=\langle
e_{i_1}\dots e_{i_s}v, \; i_1+\dots+i_s=n\rangle.
$$
$V_n(h,c)$ is an eigenspace of $e_0$ with the eigenvalue $(h+n)$:
$$
e_0(e_{i_1}\dots e_{i_s}v)=(h+i_1+\dots+i_s)e_{i_1}\dots e_{i_s}v.
$$
Besides of this $zw=cw$ for any $w \in V(h,c)$.
\begin{definition}
A vector $w\in V(h,c)$ is called singular, if $e_iw=0$ for $i <0$.
\end{definition}
A homogeneous singular vector $w \in V_n(h,c)$ of degree $n$
generates in $V(h,c)$  a submodule isomorphic to $V(h+n,c)$.

It is well-known \cite{FeFu} that there is a singular vector $w
\in V_n(h,c)$ of degree $n$ if and only if, for some positive
integers $p$ and $q$ and complex number $t$, we have $n=pq$ and
\begin{equation}
\begin{split}
c=c(t)=13+6t+6t^{-1},\\
h=h_{p,q}(t)=-\frac{p^2-1}{4}t-\frac{pq-1}{2}-\frac{q^2-1}{4}t^{-1}.
\end{split}
\end{equation}
The singular vector $w_{p,q}(t)$ of degree $n=pq$ is unique up to
scalar multiplication. One can write $w_{p,q}(t)$ as a continuous
function of $t$ \cite{Fu2}:
$$
w_{p,q}(t)=S_{p,q}(t)v=\sum_{|I|=pq}a_I^{p,q}(t)e_Iv=
\sum_{i_1+\dots+i_s=pq}a_{i_1,\dots,i_s}^{p,q}(t)e_{i_1}\dots
e_{i_s}v.
$$
The coefficients $a_I^{p,q}(t)$ depend polynomially on $t$ and
$t^{-1}$. We assume that the coefficient $a_{1,\dots,1}^{p,q}(t)$
is equal to $1$.
It is natural to consider the last sum over
ordered partitions $i_1\ge i_2\ge \dots \ge i_s \ge 1$, for
instance

\begin{equation}
\begin{split}
S_{1,1}(t)=e_1,\;\;S_{2,1}(t)=e_1^2+te_2,\;\;S_{3,1}(t)=e_1^3+4te_2e_1+(4t^2+2t)e_3,\\
S_{4,1}(t)=e_1^4+10te_2e_1^2+6t^2e_2^2+(24t^2+10t)e_3e_1+(36t^3+24t^2+6t)e_4.
\end{split}
\end{equation}
It is still unclear how to write out the general formula for all
$S_{p,q}(t)$ with the ordering $i_1\ge i_2\ge \dots \ge i_s \ge
1$.

\section{$L_1$-resolution and cohomology of thread modules}
Let us consider the Verma module $V(0,0)$. We fix the value
$t=-\frac{3}{2}$ and denote $S_{p,q}\left({-}\frac{3}{2}\right)$
by $S_{p,q}$ for simplicity.

\begin{proposition}[Kac \cite{Kac}, Feigin-Fuchs \cite{FeFu, FeFu2}]
The module $V(0,0)$ has a singular vector $w_n$ of degree $n$ (at
the level $n$) if and only if $n$ is a pentagonal number
$n=\frac{3k^2\pm k}{2}$.
\end{proposition}

Denote by $V\left(\frac{3k^2\pm k}{2}\right)$ the submodule of
$V(0,0)$ generated by the singular vector of degree $\frac{3k^2\pm
k}{2}$.

\begin{proposition}
We have the following properties:

1) the sum $V(1)+V(2)$ is the subspace of codimension one in
$V(0,0)$;

2) $V\left(\frac{3k^2- k}{2}\right)\cap V\left(\frac{3k^2+
k}{2}\right)=V\left(\frac{3(k{+}1)^2{-} (k{+}1)}{2}\right)+
V\left(\frac{3(k{+}1)^2{+} (k{+}1)}{2}\right),k {\ge} 1$.
\end{proposition}

\begin{example}
One can verify the following equality in the universal enveloping
algebra $U(L_1)$ that illustrates the property considered above:
$$
\left(e_1^3{-}6e_2e_1{+}6e_3\right)\left(e_1^2{-}\frac{2}{3}e_2\right)=
\left(e_1^4{-}\frac{20}{3}e_2e_1^2{+}4e_2^2{+}4e_3e_1{-}4e_1\right)e_1.
$$
It means that
$$
S_{3,1}S_{1,2}=S_{1,4}S_{1,1}.
$$
\end{example}

The inclusions of submodules $V\left(\frac{3k^2\pm k}{2}\right)$
considered above provides  us with the exact sequence of
Vir-modules \cite{RochWall, FeFu, FeFu2}
\begin{equation}
\label{resolution}
\begin{split}
\begin{CD}
{\dots}{\rightarrow} V(\frac{3(k{+}1)^2{-}
(k{+}1)}{2},0){\oplus}V(\frac{3(k{+}1)^2{+} (k{+}1)}{2},0)
@>{\delta_{k{+}1}}>> V(\frac{3k^2{-}
k}{2},0){\oplus}V(\frac{3k^2{+} k}{2},0){\rightarrow}{\dots}\end{CD}\\
\begin{CD}
\dots @>{\delta_3}>> V(5,0){\oplus}V(7,0) @>{\delta_2}>>
V(1,0){\oplus}V(2,0) @>{\delta_1}>> V(0,0) \rightarrow {\mathbb C}
\rightarrow 0
\end{CD},
\end{split}
\end{equation}
where $\delta_k$ are defined with the help of operators $S_{p,q}
\in U(L_1)$:
\begin{equation}
\begin{split}
\delta_{k+1}\left( \begin{array}{c}x \\y
\end{array}\right)=\left(\begin{array}{cc}S_{1,3k{+}1} & S_{2k{+}1,2}\\
{-}S_{2k+1,1}& {-}S_{1,3k{+}2} \end{array}\right)\left(
\begin{array}{c}x \\y \end{array}
\right), \; k \ge 1;\\
\delta_1\left(
\begin{array}{c}x \\y \end{array}\right)=\left(S_{1,1}, S_{1,2}\right)\left(
\begin{array}{c}x \\y \end{array}\right),
\end{split}
\end{equation}
$\varepsilon$ is the canonical surjective homomorphism.
\begin{theorem}[Rocha-Carridi-Wallach \cite{RochWall}, Feigin-Fuchs \cite{FeFu}]
The exact sequence (\ref{resolution}) regarded as the sequence of
$L_1$-modules is a free resolution of the one-dimensional trivial
$L_1$-module ${\mathbb C}$.
\end{theorem}
\begin{corollary}
Let $M$ be a $L_1$-module. Then the cohomology $H^*(L_1,M)$ is
isomorphic to the cohomology of the following complex
\begin{equation}
\label{resolution2}
\begin{CD}
{\dots}@<{d_{k{+}1}}<< M{\oplus}M @<{d_k}<< M{\oplus}M
@<{d_{k{-}1}}<< {\dots} @<{d_1}<< M{\oplus}M @<{d_0}<< M
\end{CD},
\end{equation}
with differentials
\begin{equation}
\label{differential}
\begin{split}
d_k\left(\begin{array}{c}m_1 \\m_2 \end{array}\right)=\left(\begin{array}{cc}S_{1,3k{+}1} & {-}S_{2k+1,1}\\
S_{2k{+}1,2}& {-}S_{1,3k{+}2} \end{array}\right)\left(\begin{array}{c}m_1 \\m_2 \end{array}\right), \; k \ge 1;\\
d_0(m)=\left(\begin{array}{c}S_{1,1}m \\S_{1,2}m
\end{array}\right),\;\; m, m_1, m_2 \in M.
\end{split}
\end{equation}
\end{corollary}
\begin{definition}[\cite{FeFuRe}]
An infinite dimensional thread $L_1$-module is a ${\mathbb
Z}$-graded $L_1$-module $V=\oplus_{j\in {\mathbb Z}} V_j$ such as
$$
\dim V_j=1,\quad e_i V_j \subset V_{i+j}, \; \forall e_i {\in}
L_1, j \in {\mathbb Z}.
$$
\end{definition}

One can fixe an infinite basis $f_j, \;j \in {\mathbb Z},$ such
that $f_j$ spans  $V_j$.

\begin{example}
Let us define a thread $L_1$-module $A_{\alpha}$ by its basis
$f_i,\; i \in {\mathbb Z},$ and some parameter $\alpha \in
{\mathbb K}$:
$$
e_1f_j=f_{j+1}, \quad e_2f_j=\alpha f_{j+2},\quad \forall j \in
{\mathbb Z}.
$$
It is evident that $e_if_j=0, i\ge 3, \forall j$.
\end{example}

\begin{example}
Another thread  $L_1$-module $F_{\lambda,\mu}$ came from the
well-known infinite dimensional representation of the Witt algebra
in the tensor densities \cite{Fu}:
$$ e_if_j=\left(j+\mu-\lambda(i+1) \right) f_{i+j}, \forall i \in
{\mathbb N}, j \in {\mathbb Z},
 $$
where $\lambda, \mu \in {\mathbb K}$ are two parameters.
\end{example}
For a given $L_1$-module $M$ the elements $S_{p,q} \in U(L_1)$ are
well-defined operators acting on $M$. For instance if
$M=F_{\lambda,\mu}$ then
$$
S_{2,1}f_j=\left((j{+}\mu{-}2\lambda)
(j{+}1{+}\mu{-}2\lambda){-}\frac{3}{2}(j{+}\mu{-}3\lambda)\right)f_{j{+}2}.$$
Let us introduce the numbers $\sigma_{p,q}(j) \in {\mathbb K}$
such that
$$
S_{p,q}f_j=\sigma_{p,q}(j)f_{j{+}pq}.
$$
\begin{corollary}[Feigin, Fuchs, \cite{FeFu}]
Let $M=\oplus_i M_i$ be a thread $L_1$-module over a field
${\mathbb K}$. Then the homogeneous cohomology $H^*_s(L_1,M)$ is
isomorphic to the cohomology of the following complex:
\begin{equation}
\label{resolution3}
\begin{CD}
{\dots}@<{D_{k{+}1}}<< {\mathbb K}{\oplus}{\mathbb K} @<{D_k}<<
{\mathbb K}{\oplus}{\mathbb K}@<{D_{k{-}1}}<< {\dots}  @<{D_1}<<
{\mathbb K}{\oplus}{\mathbb K} @<{D_0}<< {\mathbb K}
\end{CD},
\end{equation}
were differentials $D_k$ are defined by the matrices
\begin{equation}
D_k{=}\left(\begin{array}{cc}\sigma_{1,3k{+}1}(s{+}e({-}k)) & {-}\sigma_{2k+1,1}(s{+}e({-}k)\\
\sigma_{2k{+}1,2}(s+e(k))& {-}\sigma_{1,3k{+}2}(s+e(k))
\end{array}\right), k {\ge} 1;
D_0{=}\left(\begin{array}{c}\sigma_{1,1}(s) \\\sigma_{1,2}(s)
\end{array}\right).
\end{equation}
\end{corollary}
\begin{proof}
The homogeneous cohomology $H^*_s(L_1,M)$ is isomorphic to the
cohomology of the following subcomplex of (\ref{resolution2}).
\begin{equation}
\begin{split}
\begin{CD}
{\dots}@<{d_{k{+}1}}<< M_{s{+}e({-}k{-}1)}{\oplus}M_{s{+}e(k{+}1)}
@<{d_k}<< M_{s+e({-}k)}{\oplus}M_{s+e(k)}@<{d_{k{-}1}}<<
{\dots}\end{CD}\\
\begin{CD}
{\dots} @<{d_2}<< M_{s+5}{\oplus}M_{s+7} @<{d_1}<<
M_{s+1}{\oplus}M_{s+2} @<{d_0}<< M_s
\end{CD},
\end{split}
\end{equation}
with differentials (\ref{differential}) $d_k$ restricted on
$M_{s+e({-}k)}{\oplus}M_{s+e(k)}$. In its turn each subspace
$M_{s+e(\pm k)}$ is isomorphic to ${\mathbb K}$.
\end{proof}
At the time of writing the article \cite{FeFu} Feigin and Fuchs
were not able to write the general formula for all operators
$S_{p,q} \in U(L_1)$. However they managed to find explicit
expressions for all entries $\sigma_{p,q}(j) \in {\mathbb K}$ for
both examples $F_{\lambda, \mu}$ and $A_{\alpha}$ of thread
modules considered above. In particular for the $L_1$-module
$A_{1,\alpha}$ the matrix $D_k$ will have the following form:
$$
D_k{=}\left(\begin{array}{cc} \sigma_{1,3k{+}1} & {-}\sigma_{2k+1,1}\\
\sigma_{2k{+}1,2}& {-}\sigma_{1,3k{+}2}
\end{array}\right),
$$
where in particular
\begin{equation}
\sigma_{2k{+}1,1}=\prod_{i=1}^{k}\left(1-6\alpha i^2 \right).
\end{equation}
The corresponding formulas for the module $F_{\lambda, \mu}$ were
also obtained by Feigin and Fuchs \cite{FeFu}.

It is useful to consider also finite-dimensional thread
$L_1$-modules of the following type (short thread modules in
\cite{FeFuRe})
$$
V^{m,n}=\oplus_i V_i, \quad V_i=0,\; i < m \; {\rm or}\; i >
n,\quad \dim V_i =1,\; m \le i \le n,
$$
for some integers $m,n$.

For a given infinite dimensional thread module $V=\oplus_i V_i$
one can consider its so-called subquotient $V^{m,n}$:
$$
V^{m,n}=\left(\oplus_{i \ge m} V_i\right)/\left(\oplus_{i \ge n}
V_i\right)
$$
that is a finite-dimensional thread module: $\dim
V^{m,n}=n{-}m{+}1$.
\begin{example}
The subquotient  $A^{m,n}_{\alpha}$ can be defined by its basis
$f_m,\dots,f_n$ and non-trivial structure relations
$$
e_1f_j=f_{j+1},\; j=m,\dots,n{-}1,\; e_2f_j=\alpha f_{j+2},\;
j=m,\dots,n{-}2.
$$
\end{example}

\section{Special thread $L_1$-module $\tilde M$}
\label{new_module}
Let us consider a new thread $L_1$-module
$\tilde M$ defined by its infinite basis $\{f_j, j \in {\mathbb
Z}\}$.

\begin{equation}
\label{new_module}
e_if_j= \left\{\begin{array}{lll}
   jf_{i+j}, & j \ge 0; &\\
   (i+j)f_{i+j}, & i+j \le 0,&j<0;  \\
   f_{i+j},& i+j>0,& j<0.\\
   \end{array} \right .
\end{equation}
The verification that the formulas (\ref{new_module}) define a
$L_1$-module is straight forward. One can think of $\tilde M$ as
the result of "gluing together" of two modules: the quotient of
$F_{{-}1,1}$ with submodule of $F_{0,0}$.
\begin{remark}
The module $\tilde M=\oplus_{i} \tilde M_i$ is decomposable as the
direct sum of two $L_1$-modules
\begin{equation}
\label{decomposition} \tilde M=\tilde M_0 \oplus \left(
\oplus_{i\ne 0} \tilde M_i \right),
\end{equation}
where $\tilde M_0$ is trivial one-dimensional $L_1$-module and
$\tilde M_{\ne 0}= \oplus_{i\ne 0}\tilde M_i$ is infinite
dimensional cyclic, i.e.
\begin{equation}
\label{cyclic}
e_1\tilde M_i=\tilde M_{i+1}, i \le -2, \;
e_2\tilde M_{{-}1}=\tilde M_1, \; e_1\tilde M_i=\tilde M_{i+1}, i
\ge 1.
\end{equation}
\end{remark}
Let us consider a subquotient $\tilde M_{\ne 0}^{m,n}={\rm
Span}\;(f_m,\dots,f_{{-}1},f_1,\dots,f_n)$ that can be defined by
$$
e_if_j= \left\{\begin{array}{lll}
   jf_{i+j}, & j > 0, & i{+}j\le n;\\
   (i+j)f_{i+j}, & n \ge i+j > 0,& m \le j<0;  \\
   f_{i+j},& n \ge i+j>0,& m \le j<0.\\
   \end{array} \right.
$$
The module $\tilde M_{\ne 0}^{m,n}$ defines a representation $L_1
\to LT_{n-m+1}$ and hence a formal connection $A_{\tilde M}$. We
present below an example of formal connection $A_{\tilde M}$ with
$m={-}2, n=3$.
$$
\label{new_connect} A_{\tilde M}=\left(\begin{array}{cccccccc}
  0      & 0 & 0 & 0 & 0   & 0 &0  \\
  {-}2e^1 & 0 & 0 & 0 & 0   & 0  &0  \\
  {-}e^2  & {-}e^1  & 0 & 0 & 0      &0 &0   \\
  e^4      &  e^3      & e^2      & 0 & 0          & 0 &0   \\
  e^5      & e^4      & e^3      & e^1 & 0 & 0 & 0 \\
  e^6 & e^5      & e^4      & e^2      & 2e^1 & 0          & 0  \\
  e^7 & e^6      & e^5 & e^3      &  2e^2 & 3e^1  &  0 \\
  \end{array}\right),
$$
Now we want to establish the uniqueness in some sense of the
module $\tilde M$: up to an isomorphism $\tilde M$ is a unique
thread $L_1$-module with the decomposition (\ref{decomposition}),
where the first summand is one-dimensional trivial module and the
second one is cyclic in the sense (\ref{cyclic}).

It will be convenient for us to consider a new basis of $L_1$:
$$
\label{ind_defin} \tilde e_1 =e_1, \quad \tilde e_i =
6(i{-}2)!e_i.
$$
Now we have in particular that
\begin{equation}
[\tilde e_1,\tilde e_i]=\tilde e_{i{+}1}, \; i \ge 2.
\end{equation}
It was proved by Benoist \cite{Benoist} that $L_1$ is generated by
two elements $\tilde e_1, \tilde e_2$ with the following two
relations on them
\begin{equation}
\label{def_relat} [\tilde e_2, \tilde e_3]=\tilde e_5, \; [\tilde
e_2,\tilde e_5]=\frac{9}{10} \tilde e_7,
\end{equation}
where $\tilde e_3, \tilde e_5, \tilde e_7$ are defined by
(\ref{ind_defin}).

 Hence the defining relations (\ref{def_relat}) will give us
the following set of equations on a $L_1$-module $V={\rm Span}\;
(f_m,\dots,f_n)$:
\begin{eqnarray*}
R^5_i:\quad\quad ([\tilde e_2, \tilde e_3]-\tilde e_5)f_i=0, \;\; i=m,\dots,n{-}5,\\
R^7_j:\quad ([\tilde e_2, \tilde e_5]-\frac{9}{10}\tilde
e_7)f_i=0, \; i=m,\dots,n{-}7.
\end{eqnarray*}
Obviously if $e_1f_i \ne 0$ in a thread module $M$ one can
consider the vector $f_{i{+}1}'=e_1f_i$ instead of $f_{i{+}1}$.
\begin{theorem}
Let $M^{m,n}$ be a $(n{-}m{+}1)$-dimensional thread $L_1$-module
defined by its basis $f_i, i=m,\dots,{-}1, 0,1,\dots, n,$ with
$n{-}m{+}1 \ge 11$ such that:
\begin{eqnarray*}
\tilde e_1f_i=f_{i+1}, \;i=m,\dots,{-}2,1,\dots,n{-}1;\\
\tilde e_1f_{{-}1}=\tilde e_1f_0=0, \; \tilde e_2f_j=b_jf_{j+2},
\;
j=m,\dots,n{-}2, \quad b_{{-}2}=0, \; b_{{-}1} \ne 0.\\
\end{eqnarray*}
Then $M^{m,n}$ is isomorphic to the module $\tilde M^{m,n}$.
\end{theorem}
\begin{proof}

2) $R^5_{{-}2},R^5_{{-}1},R^5_{0},R^7_{{-}2}$ will have the
following form:
\begin{equation}
\label{dva_nulya}
\begin{split}
{-}b_{{-}1} b_{1}{-} b_{{-}2} b_{0}{+}3b_{{-}1}=0,\\
2 b_{{-}1} b_{2}{-} b_{{-}1} b_{1}{-} b_{{-}1}=0,\\
2 b_{0} b_{3}{-} b_{0} b_{2}{-} b_{0}=0,\\
{-}3 b_{{-}1} b_{3}{-}b_{{-}2} b_{0}{+}\frac{9}{2} b_{{-}1}=0.\\
\end{split}
\end{equation}
Now after rescaling the basic vectors if necessary one can assume
that
$$
b_{{-}2}=0,\;b_{{-}1}= b_{0} =1.
$$

The system (\ref{dva_nulya}) has the unique solution
$$
b_{1}=3, b_{2}=2, b_{3}=\frac{3}{2}.
$$

In order to find $b_{{-}3}$ we have to consider the following
equations:
\begin{equation}
\begin{split}
b_{{-}2} b_{0}{-} b_{{-}3} b_{{-}1}{-}3b_{{-}1}=0,\\
3b_{{-}1}b_{2}{-}b_{{-}3}b_{{-}1}{-}9b_{{-}1}=0.\\
\end{split}
\end{equation}
Evidently we have the answer $b_{{-}3}={-}3$.

For $b_{{-}4}$ we have two new additional equations:
\begin{equation}
\begin{split}
2b_{{-}4} b_{{-}1}{-} b_{{-}3} b_{{-}1}{+}b_{{-}1}=0,\\
{-}b_{{-}1}b_{1}{+}3b_{{-}4}b_{{-}1}{+}9b_{{-}1}=0.\\
\end{split}
\end{equation}
It follows that $b_{{-}4}={-}2$.

Again one can remark that there are two equations on $b_{{-}5}$:
\begin{equation}
\begin{split}
2b_{{-}5} b_{{-}2}{-} b_{{-}4} b_{{-}2}{+}b_{{-}2}=0,\\
{-}b_{0}b_{{-}2}{-}3b_{{-}5}b_{{-}1}{-}\frac{9}{2}b_{{-}1}=0.\\
\end{split}
\end{equation}
And we have $b_{{-}5}={-}\frac{3}{2}$.

Proceeding in the same way we came to the thread $L_1$-module:
$$
\begin{tabular}{|c|c|c|c|c|c|c|c|c|c|c|}
\hline
&&&&&&&&&&\\[-10pt]
 $b_m$ & $\dots$ & $b_{{-}4}$ &$b_{{-}3}$ & $b_{{-}2}$ & $b_{{-}1}$
&
 $b_{0}$ & $b_{1}$ & $b_2$ & $\dots$ & $b_{n{-}2}$ \\
&&&&&&&&&&\\[-10pt]
\hline
&&&&&&&&&&\\[-10pt]
 ${-}6/(m{+}1)$ &\dots & ${-}2$ & ${-}3$& $0$& $1$ & $1$ &
$3$ &$2$&$\dots$&$6/(n{-}1)$ \\
&&&&&&&&&&\\
\hline
\end{tabular}
$$
One can verify that this module is isomorphic to $\tilde M^{m,n}$.
\end{proof}
\begin{corollary}
\label{rigidity} Let $A'$ be an arbitrary formal connection
corresponding to the trivial Massey product $\langle
\underbrace{e^1,\dots,e^1}_m,e^2,
\underbrace{e^1,\dots,e^1}_n,\rangle, \; m+n=2k-1$. Then for an
arbitrary entry $a'(i,j), i < j$ of $A'$ we have the following
property
$$
a'(i,j)=a_{\tilde M^{m,n}_{\ne 0}}(i,j)+\delta(i,j),
$$
where $\delta(i,j)$ is a linear combination of $e^r$ with $r <
\deg a_{\tilde M^{m,n}_{\ne 0}}(i,j)$ and $a_{\tilde M^{m,n}_{\ne
0}}(i,j)$ is an entry of the formal connection $A_{\tilde
M^{m,n}_{\ne 0}}$ that corresponds to the module $\tilde M_{i\ne
0}^{m,n}$.
\end{corollary}
\begin{remark}
We will consider in the sequel both modules $\tilde M^{m,n}$ and
$\tilde M_{i\ne 0}^{m,n}$. The explanation is very simple: we need
the module $\tilde M_{i\ne 0}^{m,n}$ for the construction of
Massey products but in the same time it is more convenient to
compute the cohomology $H^*(L_1, \tilde M^{m,n})$ because of the
very simple graded thread structure of $\tilde M^{m,n}$. The
relation between $H^*(L_1, \tilde M^{m,n})$ and $H^*(L_1, \tilde
M^{m,n}_{\ne 0})$ is very simple:
$$
H^*(L_1, \tilde M^{m,n})=H^*(L_1, \tilde M^{m,n}_{\ne 0})\oplus
\tilde M_0 \otimes H^*(L_1).
$$
\end{remark}

\section{Benoit-Saint-Aubin formula and thread module $\tilde M$.}

Now it is clear that in order to have  the general combinatorial
formula for singular Virasoro vectors in terms of operators
$S_{p,q}$ it is useful to consider the sums over all unordered
sequences $pq=i_1+\dots+i_s$ by positive integers $i_1,\dots,i_s$.

The first important step in this direction was made by Benoit and
Saint-Aubin.
\begin{theorem}[Benoit, Saint-Aubin \cite{BenSA}]
\begin{equation}
\label{Benoit_SA}
\begin{split}
S_{p,1}(t)=\sum_{\begin{array}{c}i_1,\dots,i_s\\i_1+\dots+i_s=p
\end{array}.}c_p(i_1,\dots,i_s)t^{p-s}e_{i_1}\dots
e_{i_s},\\
S_{1,q}(t)=\sum_{\begin{array}{c}i_1,\dots,i_s\\i_1+\dots+i_s=q
\end{array}.}c_q(i_1,\dots,i_s)t^{-q+s}e_{i_1}\dots
e_{i_s},\\
\end{split}
\end{equation}
where the sums are all over all partitions of $p$ and $q$ by
positive numbers without any ordering restriction, and the
coefficients are defined by
\begin{equation}
c_r(i_1,\dots,i_s)=\prod_{\begin{array}{c}1\le k < r\\ k\ne
i_1+\dots+i_l, \;l=1,\dots,s{-}1.
\end{array}.} k(r-k)
\end{equation}
\end{theorem}
\begin{example}
\begin{equation}
\begin{split}
S_{1,1}(t)=e_1,\;\;S_{2,1}(t)=e_1^2+te_2,\;\;S_{3,1}(t)=e_1^3+t(2e_1e_2+2e_2e_1)+4t^2e_3,\\
S_{4,1}(t)=e_1^4+t(3e_1^2e_2+4e_1e_2e_1+3e_2e_1^2)+t^2(12e_1e_3+9e_2^2+12e_3e_1)+36t^3e_4.
\end{split}
\end{equation}
\end{example}
Later an important results on the structure result on the
structure of singular Virasoro vectors s were obtained by Bauer,
Di Francesco, Itzykson and Zuber \cite{BdFIZ} and others,
especially interesting is interpretation of singular vectors in
terms of Jack symmetric polynomials \cite{MimYam}.

 Now we are going to study the action of operators
$S_{p,1}(t)$ on the infinite dimensional thread module $\tilde M$.
\begin{example}
$$S_{3,1}(t)f_{-2}=\left(e_1^3+t(2e_1e_2+2e_2e_1)+4t^2e_3\right)f_{-2}.$$
Obviously $e_1^3f_{-2}=e_2f_{-2}=0$. Hence
$$
S_{3,1}(t)f_{-2}=(-2t+4t^2)f_1=4t\left(t-\frac{1}{2}\right)f_1
=(2!)^2\prod_{i=1}^2\left(t+\frac{i-2}{i(3-i)}\right)f_1.
$$
\end{example}
\begin{lemma}
Let $j<0$ and $p+j>0$ then
$$
S_{p,1}(t)f_j=(p-1)!^2\prod_{i=1}^{p-1}\left(t+\frac{i+j}{i(p-i)}\right)f_{j+p}.
$$
\end{lemma}
\begin{proof}
\begin{proposition}Let $j<0$ and $i_1+\dots +i_s+j>0, i_l>0$. Then
we have the following formula
$$
e_{i_1}\dots e_{i_{s-1}}e_{i_s}f_j=\prod_{l=2}^{s}(i_l+\dots
+i_s+j)f_{j+i_1+\dots +i_s}
$$
\end{proposition}
\begin{proof}
There are three possibilities:

1) $i_s+j < 0$ and it follows that  $\exists m, s\ge m \ge 2$ such
that
$$
i_s+\dots+i_m+j < 0, \quad  {\rm and} \quad i_s+\dots+i_{m-1}+j
\ge 0.
$$
Hence one can see that
\begin{equation}
\begin{split}
e_{i_{m}}\dots e_{i_s}f_j=(i_s+j)\dots(i_s+\dots
+i_{m}+j)f_{j+i_s+\dots +i_{m}}, \\
e_{i_{m-1}}f_{j+i_s+\dots +i_{m}}=f_{j+i_s+\dots +i_{m-1}},\\
e_{i_1}\dots e_{i_{m-1}}f_{j+i_1+\dots +i_{m}}=(i_s+\dots
+i_{2}+j)\dots(i_s+\dots +i_{m-1}+j)f_{j+i_s+\dots +i_1}.
\end{split}
\end{equation}

2) $i_s+j >  0$ then we have
\begin{equation}
\begin{split}
e_{i_s}f_j=f_{j+i_s}, \\
e_{i_1}\dots e_{i_{s-1}}f_{j+i_s}=(i_s+\dots
+i_{2}+j)\dots(i_s+j)f_{j+i_s+\dots +i_1}.
\end{split}
\end{equation}

3) $i_s+j =  0$ then we have $e_{i_s}f_j=(i_s+j)f_{j+i_s}=0$.
\end{proof}
Let us introduce new positive integers
$$
k_1=i_s, \; k_2=i_s+i_{s-1}>k_1, \;\dots,\; k_{s-1}=i_s+\dots+i_2
>k_{s-2}>\dots > k_1.
$$
Also one can remark that
$$
c_p(i_1,\dots,i_s)=\frac{(p-1)!^2}{i_1(i_1{+}i_2){\dots}
(i_1{+}{\dots}{+}i_{s{-}1})(p{-}i_1)(p{-}i_1{-}i_2){\dots}(p{-}i_1{-}{\dots}{-}i_{s{-}1})}.
$$
Now one can consider the Benoit-Saint-Aubin formula
(\ref{Benoit_SA}). It follows that we have the formula for
$S_{p,1}(t)f_j$:
\begin{equation}
\begin{split}
S_{p,1}(t)f_j=\\=\sum_{\begin{array}{c}(i_1{,}{\dots}{,}i_s){\in}
{\mathbb N}^s,\\i_1{+}{\dots}{+}i_s{=}p\end{array}}
\frac{t^{p-s}(p{-}1)!^2(i_s{+}j){\dots}(i_s{+}{\dots}{+}i_2{+}j)}{i_1{\dots}
(i_1{+}{\dots}{+}i_{s{-}1})(p{-}i_1){\dots}(p{-}i_1{-}{\dots}{-}i_{s{-}1})}f_{j+p}=\\
=(p{-}1)!^2\sum_{1{\le} k_1{<}\dots{<}k_{s{-}1}{\le}p{-}1}
t^{p-s}\frac{(k_1{+}j)}{k_1(p{-}k_1)}{\dots}\frac{(k_{s{-}1}{+}j)}{k_{s{-}1}(p{-}k_{s{-}1})}f_{j+p}=\\
=(p{-}1)!^2\prod_{i=1}^{p-1}\left(
t+\frac{i{+}j}{i(p{-}i)}\right)f_{j+p}.
\end{split}
\end{equation}
\end{proof}
Let us introduce a polynomial $F_{j,p}(t)$:
$$
F_{j,p}(t)=(p{-}1)!^2\prod_{i=1}^{p-1}\left(t+\frac{i+j}{i(p-i)}\right).
$$
\begin{proposition}
Let $j<0$ and $p \ge 2$, then
\begin{equation}
\begin{split}
F_{j,p}\left(-\frac{3}{2}\right) \ne 0; \quad
F_{j,p}\left(-\frac{2}{3}\right)=0  \Leftrightarrow p+3j=1.
\end{split}
\end{equation}
\end{proposition}

\begin{proof}
Let us suppose that $F_{j,p}\left(-\frac{3}{2}\right)=0$. It means
that
$$3i^2+(2{-}3p)i+2j=0$$
for some positive integer $i, 1\le i\le p{-}1$. For the roots
$i_1, i_2$ of this square equation we have $3i_1i_2=2j < 0$. On
the another hand $i_1+i_2=\frac{3p-2}{3}$ and it implies that for
the positive root $i_2$ we have $i_2 > \frac{3p-2}{3}$ that is
impossible in the virtue of $i_2 \le p-1$.

Proceeding in the analogous way in the case
$F_{j,p}\left(-\frac{2}{3}\right)=0$ one can easily see that for
the positive root $i_2$ of the corresponding square equation
$2i^2+(3{-}2p)i+3j=0$ we have $i_2 > \frac{2p-3}{2}$. It follows
that $i_2=p{-}1$ and $p{-}1{+}3j=0$.
\end{proof}

\section{Main theorem}
\begin{theorem}
\label{main_th}
The cohomology $H^*(L_1)$ is generated by two
elements $e^1, e^2 \in H^1(L_1)$ by means of two series of
non-trivial Massey products. More precisely the recurrent
procedure is organized as follows:

1) we take $e^1$ and $e^2$ as a basis of $H^1(L_1)$;

2) the triple Massey product $\langle e^1,e^2,e^2 \rangle$
 is single-valued and determines non-trivial cohomology
class $g_-^{2}{=}\langle e^1,e^2,e^2 \rangle \in H^{2}_{5}(L_1)$;

3) the $5$-fold product $\langle e^1,e^2,e^1,e^1,e^2 \rangle$ is
non-trivial and it is an affine line $\{g_+^{2}+t g_-^{2}, t \in
{\mathbb K}\}$ on the plane  $H^{2}(L_1)$, where $g_+^{2}$ denotes
some generator from $H^{2}_{7}(L_1)$. One can take an arbitrary
element $\tilde g_+^{2} \in \langle e^1,e^2,e^1,e^1,e^2 \rangle$
as the second basic element of $H^{2}(L_1)$.

Let us suppose that we have already constructed some basis $g_-^k,
\tilde g_+^k$ of $H^k(L_1), k {\ge} 2,$ such that the cohomology
class $g_-^k$ spans the homogeneous subspace
$H^k_{\frac{3k^2-k}{2}}(L_1)$. Then

4) the $(2k+1)$-fold Massey product
$$ \langle
\underbrace{e^1,\dots,e^1}_m,e^2,
\underbrace{e^1,\dots,e^1}_n,\tilde g_+^k\rangle, \; m+n=2k-1,$$
 is single-valued and determines non-trivial cohomology
class $g_-^{k+1}$ from the subspace
$H^{k{+}1}_{\frac{1}{2}(3(k{+}1)^2{-}(k{+}1))}(L_1)$.

5) the $(3k+2)$-fold product
$$ \langle
\underbrace{e^1,\dots,e^1}_{k},e^2,\underbrace{e^1,\dots,e^1}_{2k},\tilde
g^k_+ \rangle$$
 is non-trivial and it is an affine line on the two-dimensional plane $H^{k+1}(L_1)$ parallel to the
 one-dimensional subspace
 $H^{k{+}1}_{\frac{1}{2}(3(k{+}1)^2{-}(k{+}1))}(L_1)$.
One can take an arbitrary element $\tilde g_+^{k+1}$ in $\langle
\underbrace{e^1,\dots,e^1}_{k},e^2,\underbrace{e^1,\dots,e^1}_{2k},\tilde
g^k_+ \rangle$ as the second basic element of $H^{k+1}(L_1)$.
\end{theorem}
\begin{proof}

2) for an arbitrary choice of formal connection $A$ for the
product $\langle e^2,e^2,e^1\rangle$ the corresponding cocycle
will be equal to $c(A)=-e^2 \wedge e^3+\alpha d(e^3)$ for some
constant $\alpha$. Hence $\langle e^1,e^2,e^2\rangle=-[e^2 \wedge
e^3]\ne 0$ in $H^2_5(L_1)$.

3) in the second case it is more convenient to consider the
equivalent product $\langle e^1,e^2,-e^1,-2e^1,-e^2 \rangle$
instead of $\langle e^1,e^2,e^1,e^1,e^2 \rangle$. One can take the
following formal connection $A$:
$$ A=
\left(\begin{array}{cccccc}
0 & 0 & 0  &0  &0 &0    \\
  -e^2 & 0 & 0  &0  &0 &0    \\
   2e^3   & -2e^1   & 0       &0  &0&0     \\
     -e^4{-}te^2 & -e^2 &-e^1  & 0 & 0 &0       \\
      0&  e^4 & e^3 & e^2 & 0 &0 \\
      *&   e^5 & e^4 & e^3& e^1 &0          \\
\end{array}\right).
$$
The corresponding cocycle will be $c(A)=(e^2{\wedge}
e^5{-}3e^3{\wedge} e^4)+te^2{\wedge} e^3$. On the another hand for
an arbitrary defining system $A'$ the corresponding cocycle will
have the form $c(A')=(e^2{\wedge} e^5{-}3e^3{\wedge}
e^4){+}{\dots}$, where dots stand for the summands with the second
grading strictly less than $7$.

4) Let us consider a finite dimensional thread  $L_1$-module
$\tilde M_{-m-1,n}=\langle f_{-m-1},f_m,\dots,f_n\rangle$ where
$m+n=2k-1$. In the second grading $-m-1-\frac{3k^2+k}{2}$ in the
spectral sequence we have only one differential that can be non
trivial:
$$
\begin{CD}
 g_+^k \otimes f_{-m-1} @>{d_{2k+1}}>>  g_-^{k+1}
\otimes f_n .
\end{CD}
$$
We recall that
$$
d_{2k+1}\left(g_+^k \otimes f_{-m-1}\right)=c(A_{\tilde
M^{-m-1,n}_{\ne 0}})\otimes f_n,
$$
where $c(A_{\tilde M^{-m-1,n}_{\ne 0}})$ is the cocycle that
corresponds to the formal connection  $A_{\tilde M^{-m-1,n}_{\ne
0}}$ of the $(2k+1)$-fold Massey product
$$
\langle
\underbrace{e^1,\dots,e^1}_m,e^2, \underbrace{e^1,\dots,e^1}_n,
g_+^k\rangle, \; m+n=2k-1,
$$
defined by the module $\tilde M^{-m-1,n}_{\ne 0}$.
\begin{proposition}
The differential $d_{2k+1}$ (and hence the cocycle $c(A_{\tilde
M^{-m-1,n}_{\ne 0}})$) is non trivial.
\end{proposition}
\begin{proof}
We will compute the cohomology by means of the free resolution.
Here we have also only one possibly non trivial differential:
$$
\begin{CD}
 (0,f_{-m-1})\rangle @>{D_k}>>  (f_n,0),
\quad D_k=\left( \begin{array}{cc}0& 0\\-S_{2k{+}1, 1} & 0\\
\end{array}\right).
\end{CD}
$$
It follows that
$$
S_{2k+1,1}(f_{-m-1})=F_{-m-1,2k+1}\left(-\frac{3}{2}\right)f_{n}
\ne 0.
$$
\end{proof}
\begin{corollary}
It exists a formal connection $A$ for the Massey product
$$ \langle
\underbrace{e^1,\dots,e^1}_m,e^2, \underbrace{e^1,\dots,e^1}_n,
g_+^k\rangle, \; m+n=2k-1,$$ such that $c(A)=g_-^{k+1}$.
\end{corollary}

\begin{proposition}
The product $\langle \underbrace{e^1,\dots,e^1}_m,e^2,
\underbrace{e^1,\dots,e^1}_n, g_+^k\rangle, \; m+n=2k-1,$ is non
trivial and single-valued, moreover one can replace in it the
element $g_+^k$ by $\tilde g_+^k=g_+^k+tg_-^k$ with arbitrary $t$.
\end{proposition}
\begin{proof}
It directly follows from the Corollary \ref{rigidity}: for an
arbitrary formal connection $A'$ that corresponds to our Massey
product we will have
$$
c(A')=c(A_{\tilde M})+\dots,
$$
where dots stand for summands with the gradings strictly less than
$\deg c(A_{\tilde M})$ and hence $[c(A')]=[c(A_{\tilde M})]$
because the grading $\deg c(A_{\tilde M})=e_-(k+1)$ is the minimal
where there is non trivial cohomology class from $H^{k+1}(L_1)$.
\end{proof}

5) Let us consider a finite dimensional graded  $L_1$-module
$\tilde M^{-2k-1,k}$. In the second grading
${-}2k{-}1-\frac{3k^2+k}{2}$ in the spectral sequence we have only
one possibly non trivial differential
$$
\begin{CD}
\langle g_+^k \otimes f_{-2k-1}\rangle @>{d_{2k+1}}>> \langle
g_-^{k+1} \otimes f_k \rangle.
\end{CD}
$$
We recall that
$$
d_{2k+1}\left(g_+^k \otimes f_{-2k-1}\right)=c(A_M)\otimes f_k,
$$
where $c(A_M)$ is the cocycle that corresponds to the formal
connection $A_M$ of the Massey product
$$ \langle
\underbrace{e^1,\dots,e^1}_{k},e^2,\underbrace{e^1,\dots,e^1}_{2k},
g^k_+ \rangle$$ defined by the module $M_{-2k-1,k}$.
\begin{proposition}
The differential $d_{3k+2}$ (and hence the cocycle $c(A_M)$) is
non trivial.
\end{proposition}
\begin{proof}
The computation by means of the free resolution gives us that in
this complex we have also the only one possibly non trivial
differential:
$$
\begin{CD}
\langle (0,f_{-2k-1})\rangle @>{D_k}>> \langle (f_k,0)\rangle,
\quad D_k=\left( \begin{array}{cc}0& 0\\-S_{2k{+}1, 1} & -S_{1, 3k+2}\\
\end{array}\right).
\end{CD}
$$
It is obvious that
$$
S_{2k+1,1}(t)f_{-2k-1}=0.
$$
On the another hand
$$
S_{1, 3k+2}\left(-\frac{3}{2} \right)f_{-2k-1}=S_{3k+2,
1}\left(-\frac{2}{3} \right)f_{-2k-1} \ne 0.
$$
\end{proof}
\begin{corollary}
It exists a formal connection $A$ for the Massey product
$$ \langle
\underbrace{e^1,\dots,e^1}_{k},e^2,\underbrace{e^1,\dots,e^1}_{2k},
g^k_+ \rangle$$such that $c(A)=g_+^{k+1}$.
\end{corollary}
\begin{proposition}
The product $\langle \underbrace{e^1,\dots,e^1}_k,e^2,
\underbrace{e^1,\dots,e^1}_{2k}, g_+^k\rangle$ is non trivial, but
not single-valued. It value will not change if we replace in it
the element $g_+^k$ by $\tilde g_+^k=g_+^k+tg_-^k$ with arbitrary
constant $t$.
\end{proposition}
\begin{proof}
It follows from the Corollary \ref{rigidity}.
\end{proof}

\end{proof}

\end{document}